\documentclass[letterpaper,11pt]{article}
\usepackage{Macros}
\usepackage{hyperref}
\title{Untilting Line Bundles on Perfectoid Spaces}
\author{Gabriel Dorfsman-Hopkins\footnote{University of Caliornia, Berkeley.  Department of Mathematics.  895 Evans Hall, Berkeley, CA 94720.
\newline \textit{Email:} gabrieldh@berkeley.edu}}
\begin{document}
\date{}
\maketitle
\begin{abstract}
  Let $X$ be a perfectoid space with tilt $X^\flat$.  We build a natural map $\theta:\Pic X^\flat\to\lim\Pic X$ where the (inverse) limit is taken over the $p$-power map, and show that $\theta$ is an isomorphism if $R = \Gamma(X,\sO_X)$ is a perfectoid ring.  As a consequence we obtain a characterization of when the Picard groups of $X$ and $X^\flat$ agree in terms of the $p$-divisibility of $\Pic X$.  The main technical ingredient is the vanishing of higher derived limits of the unit group $R^*$, whence the main result follows from the Grothendieck spectral sequence.
\end{abstract}
\section{Introduction}
Scholze \cite{sh12} introduced a class of algebro-geometric objects called \textit{perfectoid spaces} that arise naturally in the context of $p$-adic geometry.  To a perfectoid space $X$ over the $p$-adic numbers, one can functorially assign a perfectoid space $X^\flat$ in characteristic $p$, called the \textit{tilt} of $X$.  Remarkably, $X$ and $X^\flat$ share many algebraic and topological properties, enumerated through various \textit{tilting equivalences}, including that they are homeomorphic and have canonically identified \'etale sites.  Perhaps surprisingly given their similarities, it is not true in general that $X$ and $X^\flat$ have isomorphic Picard groups (see \cite[Section 6]{heu21b} or Example \ref{TateCurve} below).  The goal of this paper is to describe exactly how different these two Picard groups can be.

Observe that since $X^\flat$ is perfect, $p$ acts invertibly on its Picard group.  We will show that the only obstruction to the Picard groups of $X$ and $X^\flat$ being isomorphic is the failure of $p$ to act invertibly on $\Pic X$.  To see this we construct a map
\[\theta:\Pic X^\flat\longto\lim\Pic X\]
where the (inverse) limit is taken along the $p$-power map on $\Pic X$.  A rough outline of the construction is the following: the tilting equivalence identifies the \'etale sites of $X$ and $X^\flat$ \cite[Theorem 7.9]{sh12}, so that we can view $\sO_X^*$ and $\sO_{X^\flat}^*$ as objects in the same category.  The construction of tilting then provides an identification:
\[\sO_{X^\flat}^*\cong\lim\sO_X^*,\]
where the limit is taken along the map $f\mapsto f^p$.  One can exhibit $\theta$ by taking derived global sections of this isomorphism and analyzing the Grothendieck spectral sequence.  The obstructions to $\theta$ being an isomorphism are evidently controlled by derived limits of $\sO_X^*$ and $R^* = \HH^0(X,\sO_X^*)$.  Therefore, the main technical component of this paper is the study of these derived limits.  We do this in Section \ref{limitSection}, showing the following.
\begin{Proposition}\label{limitLemmaIntro}
  If $R$ is a perfectoid ring then $\mathbf{R}^1\lim R^*=0$.
\end{Proposition}
We prove the main theorem in Section \ref{homologicalSection}, constructing $\theta$ and using the Grothendieck spectral sequence to explicitly compute the obstructions to injectivity and surjectivity in terms of derived limits of $\sO_X^*$ and $R^*$.  This allows us to leverage Proposition \ref{limitLemmaIntro} to prove the main result.
\begin{Theorem}\label{main}
  Let $X$ be a perfectoid space with tilt $X^\flat$.  There is a natural map $\Pic X^\flat\to\lim\Pic X$ which is an isomorphism if $\Gamma(X,\sO_X)$ is a perfectoid ring.
\end{Theorem}
Composing the target of $\theta$ with projection onto the first coordinate one obtains a characterization of when the Picard groups of $X$ and $X^\flat$ agree.
\begin{Corollary}\label{corToMain}
  Let $X$ be a perfectoid space with tilt $X^\flat$ and assume that $\Gamma(X,\sO_X)$ is a perfectoid ring.  There is a natural map $\Pic X^\flat\to\Pic X$ which is an isomorphism if and only if $\Pic X$ is uniquely $p$-divisible.
\end{Corollary}
In fact, the kernel and cokernel of this map can be explicitly computed in terms of the action of $p$ on $\Pic X$.  The kernel can be identified with the $p$-adic Tate module of $\Pic X$, and the image consists of the maximal $p$-divisible subgroup of $\Pic X$.

We continue by enumerating consequences and examples, including examples where the Picard groups of $X$ and $X^\flat$ agree, and examples where they don't.  For instance, a consequence of Bhatt and Scholze's theory of prismatic cohomology is that the Picard groups of a large class of perfectoid rings are uniquely $p$-divisible \cite[Corollary 9.5]{bs19}, so that we may extract the following result.
\begin{Corollary}\label{corToPrisms}
  Let $R$ be a perfectoid ring over a perfectoid field, and let $R^\flat$ be its tilt.  Then $\Pic R^\flat\cong\Pic R$.
\end{Corollary}
In the final section we apply these results to describe an example of a perfectoid cover of an abelian variety over $\bC_p$ with good reduction, $\widetilde{A}\to A$ (as in \cite{ps16} and \cite{bghsswy18}), where the natural map $\colim\Pic A\to\Pic \widetilde{A}$ is \textit{not} surjective.  This stands in stark contrast to the cases of the perfectoid cover of projective space and perfectoid covers of toric varieties where the analogous map is an isomorphism (\cite[Theorem 3.4]{DH21} and \cite[Theorem 4.1]{dhrw20} respectively). To do this we use an analog of Theorem \ref{main} to compare the N\'eron-Severi groups of $A$, $\widetilde{A}$, and $\widetilde{A}^\flat$.  This is not exclusively a mixed characteristic phenomenon, work of Heuer \cite[Section 6.2]{heu21b} constructs such an example in characteristic $p$, using a comparison between the Picard group of $\widetilde{A}$ and the Picard group of the perfection of the special fiber of $A$.  In both cases the analysis boils down the phenomenon of Picard rank potentially jumping on a special fiber.

In \cite{DH21} we constructed $\theta$ from a geometric perspective using the notion of \textit{projectivoid geometry}.  This construction required the existence of an ample line bundle on $X^\flat$ to guarantee that one can associate line bundles on $X^\flat$ to maps to a perfectoid analog of projective space using \cite[Theorem 4.5]{DH21} which in turn are in one-to-one correspondance for $X$ and $X^\flat$ by the tilting equivalence of Scholze \cite[Proposition 6.17]{sh12}, thereby giving $\theta$.  With this we gave a geometric proof that $\theta$ was injective \cite[Theorem 5.13]{DH21}, although in less generality than this paper.
\vspace{-7pt}
\subsection{Acknowledgments}
\vspace{-7pt}
The author is indebted to Bhargav Bhatt for communicating many of the main ideas of this proof over several conversations, and to Yupeng Wang for pointing out and correcting an error in the proof of Proposition \ref{pDivisible} in an earlier draft.  The author also thanks Arthur Ogus for comments greatly simplifying the exposition, and to Ben Heuer, Martin Olsson, and Peter Wear helpful conversations.  The author is partially supported by NSF grant DMS-1646385 as part of the Research Training Group in arithmetic geometry at the University of California, Berkeley.

\section{Derived Limits of Unit Groups of Perfectoid Rings}\label{limitSection}
For what follows we will denote by $\lim$ the functor that takes an abelian group (or sheaf of abelian groups) $A$ to the inverse limit of the system $(\cdots\to A\to A)$ whose transition maps are all multiplication by $p$ (or $p$-power if the group is multiplicative).  We will denote the derived functors by $\mathbf{R}^i\lim$ and the derived inverse limit by $\mathbf{R}\lim$.  The main goal of this section is to establish the vanishing of derived limits of the unit groups of perfectoid rings, in both the integral and Tate cases.  We begin with the integral case.
\begin{Definition}[Integral Perfectoid Rings {\cite[Definition 3.5]{BMS1}}]
  A topological ring $R$ is called an \textit{integral perfectoid ring} if the following conditions hold.
  \begin{enumerate}[label ={(\arabic*)}]
    \item{
    There is some $\varpi\in R$ with $\varpi^p|p$ such that $R$ has the $\varpi$-adic topology.
    }
    \item{
    Frobenius surects on $R/pR$.
    }
    \item{
    Fontaine's period map $\bA_{\inf}(R)\to R$ is surjective with principal kernel.
    }
  \end{enumerate}
\end{Definition}
\begin{Remark}\label{BMSComparison}
  If $\varpi$ is a non-zero-divisor, then by \cite[Lemma 3.10]{BMS1}, one can replace conditions (2) and (3) with the condition that Frobenius induces an isomorphism $R/\varpi R\longtoo{\sim} R/\varpi^p R$, thereby avoiding the need to invoke $\bA_{\inf}$.  We point out that, $\varpi$ is automatically a non-zero divisor if $R$ is the ring of definition for a perfectoid Tate ring (see Definition \ref{perfTateDef} below), which is our main case of interest.
\end{Remark}
We will first show that if $R$ is integral perfectoid then $R^*$ has vanishing derived limits.  The following lemma sets up the argument.
\begin{Lemma}\label{limitLemma}
  Let $(A_i)$ be an inverse system of abelian groups and let $A = \lim_i A_i$.  Suppose the following two conditions hold.
  \begin{enumerate}[label ={(\arabic*)}]
    \item{For all $i$, $\lim A_i=0$.}
    \item{For all $i$, $\bfR^1\lim A_i = 0$.}
  \end{enumerate}
  Then $\mathbf{R}^1\lim A = 0$.
\end{Lemma}
\begin{Proof}
  We derive the identification $\lim(\lim_i A_i)\cong\lim_i(\lim A_i)$ and consider the composition of functors spectral sequence which induces the following:
  \[
  \begin{tikzcd}
    0\rar&\bfR^1\lim(\lim_i A_i)\rar & \bfR^1(\lim\circ\lim_i)(A_i)\rar\ar[d,"\cong"]&\lim(\bfR^1\lim_i A_i)\\
    0\rar & \bfR^1\lim_i(\lim A_i)\rar&\bfR^1(\lim_i\circ\lim )(A_i)\rar&\lim_i\bfR^1\lim A_i.
  \end{tikzcd}
  \]
  The bottom left corner is 0 by condition (1), and the bottom right corner is 0 by condition (2).  Therefore the top left corner is 0, but this is $\bfR^1\lim A$, so we are done.
\end{Proof}
We would like to apply this in the following case.
\begin{Lemma}\label{filtration}
  Suppose $A$ is an abelian group with a complete descending filtration whose subquotients are $p$-torsion.  Then $\bfR^1\lim A=0$.
\end{Lemma}
\begin{Proof}
  We denote the filtration by $\cdots\subseteq F_2\subseteq F_1\subseteq F_0 = A$, let $Q_i = F_i/F_{i+1}$ be the subquotients, and consider the inverse system $(\cdots\to Q_i\to Q_i)$ whose transition maps are multiplication by $p-$hence 0 since $Q_i$ is $p$-torsion.  Therefore $\lim Q_i = 0$ and the system has stable images (which are all 0), so that it satisfies the Mittag-Leffler condition and $\bfR^1\lim Q_i=0$ as well.

  Let $A_i = A/F_i$.  We have the following diagram of exact sequences for all $i$.
  \[
  \begin{tikzcd}
    0\rar &F_{i+1}\ar[d,equals]\rar& F_i\rar\dar &Q_i\rar\dar & 0\\
    0\rar &F_{i+1}\rar & A\rar & A_{i+1}\rar& 0.
  \end{tikzcd}
  \]
  Therefore by the snake lemma we obtain exact sequences
  \[0\longto Q_i\longto A_{i+1}\longto A_i\longto 0,\]
  for all $i$.  Passing to the limit along multiplication by $p$ one obtains the following long exact sequence:
  \[0\to\lim Q_i\to\lim A_{i+1}\to\lim A_i\to\bfR^1\lim Q_i\to\bfR^1\lim A_{i+1}\to\bfR^1\lim A_i.\]
  Arguing inductively we may assume that $\lim A_i = 0 = \bfR^1\lim A_i$ (since $A_1 = Q_1$), so that $\lim A_{i+1} = 0 = \bfR^1\lim A_{i+1}$.  Lastly, since we are assuming that $A$ is complete with respect to the topology induced by the $F_i$, we see that the natural map $A\to\lim_i A_i$ is an isomorphism, so that by Lemma \ref{limitLemma} we may conclude that $\bfR^1\lim A=0$.
\end{Proof}
With this in hand we can take care of the main result of this section in the integral perfectoid case.
\begin{Proposition}\label{limInts}
  Let $R$ be an integral perfectoid ring.  Then $\mathbf{R}^1\lim R^{*}=0$
\end{Proposition}
\begin{Proof}
  Since $\varpi$ lives in the Jacobson radical of $R$, the projection map $R\to R/\varpi R$ remains surjective on unit groups, so that we have an exact sequence
  \[1\longto (1+\varpi R)^*\longto R^*\longto (R/\varpi R)^*\longto 1.\]
  As $R$ is perfectoid, Frobenius surjects on $R/\varpi R$ \cite[Lemma 3.9]{BMS1}, so that the $p$-power map surjects on $(R/\varpi R)^*$.  In particular, the system $(\cdots\to (R/\varpi R)^*\to(R/\varpi R)^*)$ satisfies the Mittag-Leffler condition and so $\mathbf{R}^1\lim (R/\varpi R)^* =0$. Considering the long exact sequence derived from applying $\lim$ to the sequence above therefore provides a surjection
  \[
  \begin{tikzcd}
    \mathbf{R}^1\lim(1+\varpi R)^*\ar[r,twoheadrightarrow]& \mathbf{R}^1\lim R^*,
  \end{tikzcd}
  \]
  so that it suffices to prove the source vanishes.  But $(1+\varpi R)^*$ has a complete filtration by $(1+\varpi^nR)^*$ for $n\ge0$, and the subquotients are $(1+\varpi^n R/\varpi^{n+1} R)^*$.  Since $\varpi|p$, we may observe that these subquotients are $p$-torsion.  Indeed, lifting an arbitrary element to $1+f\varpi^n\in (1+\varpi^n R)^*$, we have
  \[(1+f\varpi^n)^p = 1 + (pf\varpi^{n} + \cdots ) = 1+\varpi^{n+1}(\cdots)\equiv1\mod\varpi^{n+1}.\]
  Therefore by Lemma \ref{filtration} we deduce that $\bfR^1\lim(1+\varpi R)^*=0$ completing the proof.
\end{Proof}
For the remainder of this section we extend this result to the perfectoid Tate case.  We first recall the definition.
\begin{Definition}[Perfectoid Tate Rings]\label{perfTateDef}
  A complete topological ring $R$ is \textit{Tate} if there exists an open subring of definition $R_0\subseteq R$ whose topology is $\varpi$-adic for some unit $\varpi\in R_0\cap R^*$.  The element $\varpi$ is called a \textit{pseudouniformizer} for $R$.

  A Tate ring $R$ is called a \textit{perfectoid Tate ring} if the subring $R^\circ$ of power bounded elements in $R$ is an integral perfectoid ring.
\end{Definition}
\begin{Remark}\label{algebraicTilting}
  To a perfectoid ring $R$ one can associate a perfectoid ring $R^\flat$ of characteristic $p$ called the \textit{tilt} of $R$.  As a multiplicative monoid $R^\flat=\lim_{x\mapsto x^p} R$ so that it comes equipped with a multiplicative map $\sharp:R^\flat\to R$ which is identified with projection onto the first coordinate.  $R^\flat$ and $R$ share many algebraic and topological properties, explored through the various \textit{tilting equivalences} of Scholze \cite{sh12}.
\end{Remark}
\begin{Remark}
  By taking $\varpi = (\varpi^\flat)^\sharp$ for any pseudouniformizer $\varpi^\flat$ of $R^\flat$, one obtains a pseudouniformizer of $R$ equipped with a compatible set of $p$-power roots $\varpi^{1/p^n}$ for every $n$.  This way, the symbol $\varpi^d$ makes sense for any $d\in\bZ[1/p]$.  For what follows, we always use such a $\varpi$.
\end{Remark}
Fix a perfectoid Tate ring $R$.  We will reduce the study of the derived limits of $R^*$ to that of $R^{\circ*}$ (whose derived limits we already know vanish by Proposition \ref{limInts}) using the following exact sequence:
\begin{equation}\label{imptSequence}
  0\longto R^{\circ*}\longto R^*\longto R^*/R^{\circ*}\longto 0
\end{equation}
It suffices to show that the right hand side of this sequence is $p$-divisible.  We can in fact show something slightly stronger: that it is uniquely so.
\begin{Proposition}\label{pDivisible}
  Let $R$ be a perfectoid Tate ring.  Then $R^*/R^{\circ*}$ is uniquely $p$-divisible.
\end{Proposition}
\begin{Proof}
  We first show the $p$-divisibilty.  Fix $f\in R^*$.  We hope to show that there is some $g\in R^*$ so that $f/g^p\in R^{\circ*}$.  We first take care of the case where $f\in R^\circ$.  Let $X = \spa(R,R^\circ)$.  As $f$ is invertible, one has $|f(x)|>0$ for all $x\in X$.  Since $X$ is compact, there is some $d$ with $|\varpi^d|<|f(x)|$ for all $x\in X$.  Applying an approximation lemma \cite[Lemma 2.3.1]{cs19}, there exists some $g\in R^\flat$ so that for all $x\in X$.
  \begin{eqnarray*}
    |f(x) - g^\sharp(x)|&\le&|p|\cdot\max\{|f(x)|,|\varpi^d|\}\\
    &=&|p|\cdot|f(x)|\\
    &<&|f(x)|.
  \end{eqnarray*}
  The nonarchimedean propety of the valuations $x\in X$ therefore implies $|f(x)| = |g^\sharp(x)|$ for all $x\in X$ (cf. \cite[Remark 6.6]{sh12}).  In particular, we know $|g^\sharp(x)|>0$ for all $x\in X$, so that $g^\sharp\in R^*$ \cite[Lemma 1.4]{hu94}.  Let $u = f/g^\sharp$.  As anything in the image of $\sharp$ is a $p$-power, it suffices to show $u\in R^{\circ*}$.  Since the valuations associated to $x\in X$ are multiplicative, one has $|u(x)| = 1$ for all $x\in X$, so that $u\in R^{*}$, and $|u^{-1}(x)| = 1$ for all $x\in X$.  In particular, by \cite[Lemma 3.3(i)]{hu93}, $u,u^{-1}\in R^\circ$ so that $u\in R^{\circ*}$.  For the general $f\in R^*$, we note that there exists $d>0$ so that $\varpi^df\in R^\circ$.  By the previous paragraph there exists $g\in R^*$ so that $\varpi^df/g^p = f/(\varpi^{d/p}g)^p\in R^{\circ*}$, as desired.

  To establish uniqueness we must show that the $p$-power map is injective on $R^*/R^{\circ*}$. We first observe that the following diagram together with the snake lemma implies that $\sO_X^*/\sO_X^{+*}$ is uniquely $p$-divisible.
  \[
  \begin{tikzcd}
    &0\dar&0\dar&{}&\\
    &\mu_p\dar\ar[r,equals]&\mu_p\dar&{}&\\
    0\rar&\sO_X^{+*}\dar\rar&\sO_X^*\dar\rar&\sO_X^*/\sO_X^{+*}\rar\dar&0\\
    0\rar&\sO_X^{+*}\dar\rar&\sO_X^*\dar\rar&\sO_X^*/\sO_X^{+*}\rar&0\\
    &0&0&{}&
  \end{tikzcd}
  \]
  The long exact sequence in cohomology for the bottom 2 rows induces
  \[
  \begin{tikzcd}
    0\rar&R^*/R^{\circ *}\dar\rar&\HH^0(X,\sO_X^*/\sO_X^{+*})\rar\ar[d,"\cong"]&\HH^1(X,\sO_X^{+*})\dar\rar & \cdots\\
    0\rar&R^*/R^{\circ *}\rar&\HH^0(X,\sO_X^*/\sO_X^{+*})\rar&\HH^1(X,\sO_X^{+*})\rar & \cdots.\\
  \end{tikzcd}
  \]
  In particular, the injectivity of the vertical map on the left follows from the commutativity of the left square and the fact that the remaining three maps in the square are all injective.
\end{Proof}
Now we can prove the main result of this section in the perfectoid Tate case.
\begin{Proposition}\label{R1LimVanishes}
  Let $R$ be a perfectoid Tate ring.  Then $\mathbf{R}^1\lim R^*$ vanishes.
\end{Proposition}
\begin{Proof}
  We apply $\lim$ to Sequence (\ref{imptSequence}) to obtain the long exact sequence,
  \[\cdots\longto \mathbf{R}^1\lim R^{\circ*}\longto \mathbf{R}^1\lim R^*\longto \mathbf{R}^1\lim R^*/R^{\circ*}\longto\cdots\]
  The left term vanishes by Lemma \ref{limInts}.  Meanwhile, by Proposition \ref{pDivisible} we know $R^*/R^{\circ*}$ is $p$-divisible, so that the inverse system $(\cdots\to R^*/R^{\circ*}\to R^*/R^{\circ*})$ has surjective transition maps, and therefore satisfies the Mittag-Leffler condition.  This implies that the right term vanishes, and so the middle one must as well by exactness and we win.
\end{Proof}
\section{Untilting Line Bundles}\label{homologicalSection}
\begin{Definition}
  A perfectoid space is an adic space with an open cover by affinoid adic spaces $\spa(R,R^+)$ where $R$ is a perfectoid Tate ring, or $\spa(R,R)$ where $R$ is an integral perfectoid ring.  (In the literature the latter is often called a perfectoid formal scheme, but we will not distinguish between the two cases in what follows).
\end{Definition}
\begin{Remark}\label{GlobalTilting}
  The tilting process of Remark \ref{algebraicTilting} glues, so that to a perfectoid space $X$ one can associate its \textit{tilt}: a perfectoid space $X^\flat$ of characteristic $p$ which is locally covered by the adic spectra of the tilts of the perfectoid affinoid rings defining the charts of $X$.  The various \textit{tilting equivalences} referred to in Remark \ref{algebraicTilting} globalize as well, for example $X$ and $X^\flat$ are homeomorphic and have equivalent \'etale sites \cite{sh12}.
\end{Remark}
Let $X$ be a perfectoid space with tilt $X^\flat$.  We begin this section by constructing the desired map between the Picard groups of $X$ and $X^\flat$.  By \cite[Theorem 3.5.8]{KL16}, vector bundles on $X$ in the analytic, \'etale, pro-\'etale, and v-topologies all agree (and similarly for $X^\flat$), so there is no confusion by what we mean by the Picard group.  Denote the sheaves of units of $X$ and $X^\flat$ by $\sO_X^*$ and $\sO_{X^\flat}^*$ respectively.  As $X$ and $X^\flat$ have canonically isomorphic analytic, \'etale, pro-\'etale, and v-sites, we can view both as sheaves on $X$.  With this perspective in mind, we have an identification:
\begin{equation}\label{ID}
  \sO_{X^\tilt}^*\cong\lim\sO_X^*.
\end{equation}
The desired map is easily obtained from this identification.
\begin{Lemma}\label{thetaExists}
  There is a natural map $\theta:\Pic X^\flat\to\lim\Pic X$.
\end{Lemma}
\begin{Proof}
  Using the identification in Equation (\ref{ID}) we obtain projections $\sO^*_{X^\flat}\to\cdots\to\sO_X^*\to\sO_X^*$.  Taking the first cohomology group gives $\Pic X^\flat\to\cdots\to\Pic X\to\Pic X$, and $\theta$ arises as the universal factorization through the inverse limit.
\end{Proof}
We can explicitly study the obstructions to $\theta$ being an isomorphism by exhibiting it in a spectral sequence.  In particular, we apply derived global sections to Equation (\ref{ID}) and analyze the Grothendieck spectral sequence.  We have the following diagram of functors.
\begin{equation}\label{diagramOfFunctors}
  \begin{tikzcd}
    & \textbf{ShAb}(X)\ar[dr,"\Gamma"]&\\
    \textbf{ShAb}(X)\ar[ur,"\lim"]\ar[dr,"\Gamma"]& &\textbf{Ab}\\
    & \textbf{Ab}\ar[ur,"\lim"] &
  \end{tikzcd}.
\end{equation}
\begin{Lemma}\label{preservingInjectivity}
  $\Gamma$ takes injective sheaves to $\lim$-acyclic abelian groups, and $\lim$ takes injective sheaves to $\Gamma$-acyclic sheaves.
\end{Lemma}
\begin{Proof}
  In fact, the global sections functor takes injectives to injectives because it has an exact left adjoint (given by the constant sheaf associated to an abelian group).

  To show $\lim$ takes injectives to acyclics we fix an injective module $\sI$, and observe that for every $U$, $\sI|_U$ is injective also (since restriction to $U$ is right adjoint to extension by 0 which is exact).   The global sections of an injective sheaf is an injective abelian group and is therefore divisible.  Thus the multiplication by $p$-map on $\sI$ is surjective.  Furthermore, by injectivity we have $\HH^i(U,\sI)=0$ for all $i>0$, so that for any $i$ and $U$ the inverse system along multiplication by $p$ \[\left(\cdots\to\HH^i(U,\sI)\to\HH^i(U,\sI)\right)\]
  has surjective transition maps and therefore satisfies the ML-condition.  Therefore, applying \cite[Proposition 13.3.1]{EGA3} we observe that for all $i$ the natural map:
  \[\HH^i(X,\lim\sI)\to\lim\HH^i(X,\sI),\]
  is an isomorphism.  In particular, $\lim\sI$ is $\Gamma$-acyclic as desired.
\end{Proof}
This allows us to harness the composition of functors spectral sequence to study morphism at the center of the manuscript and prove our main theorem.
\begin{Theorem}\label{Main}
  Let $X$ be a perfectoid space with tilt $X^\flat$.  Then $\theta$ fits into an exact sequence:
  \begin{equation}\label{mainSequence}
    0\longto \bfR^1\lim\Gamma(X,\sO_X^*)\longto\Pic X^\flat\longtoo{\theta}\lim\Pic X\longto 0.
  \end{equation}
  In particular, if $\Gamma(X,\sO_X)$ is a perfectoid ring, then $\theta$ is an isomorphism.
\end{Theorem}
\begin{Proof}
  Consider the compositions in Diagram (\ref{diagramOfFunctors}).  Applying Lemma \ref{preservingInjectivity} and \cite[Tag 015M]{stacks-project} one obtains natural isomorphisms of derived functors:
  \begin{equation}\label{identifications}
    \bfR\Gamma(X,\bfR\lim(\cdot))\oo{\longleftarrow}{\sim}{}\bfR(\Gamma(X,\cdot)\circ\lim(\cdot))\cong\bfR(\lim\circ\Gamma(X,\cdot))\longtoo{\sim}\bfR\lim\bfR\Gamma(X,\cdot).
  \end{equation}
  Furthemore, the cohomology of the middle complex computes the sheaf cohomology of $\lim$ of a sheaf.  Therefore, identifying $\sO_{X^\flat}^*\cong\lim\sO_X^*$ as abelian sheaves on $X$, as well as making the appropriate identifiations using Equation (\ref{identifications}), we obtain the following two composition of functors spectral sequences:
  \begin{eqnarray*}
    E^{p,q}_2:&\HH^p(X,\mathbf{R}^q\lim\sO_X^*)&\Longrightarrow \HH^{p+q}(X,\mathbf{R}\lim\sO_X^*),\\
    E'^{p,q}_2:&\mathbf{R}^p\lim \HH^q(X,\sO_X^*)&\Longrightarrow \HH^{p+q}(X,\mathbf{R}\lim\sO_X^*).
  \end{eqnarray*}
  Considering the low degree terms gives us the following diagram whose rows are exact, and $\theta$ appears as the diagonal arrow.
  \begin{equation}\label{lowTerms}
    \begin{tikzcd}
      0\rar&\HH^1(X,\sO_{X^\tilt}^*)\ar[drr]\rar&\HH^1(X,\mathbf{R}\lim\sO_X^*)\ar[d,equal]\rar&\Gamma(X,\mathbf{R}^1\lim\sO_X^*)\rar&\HH^2(X,\sO_{X^\flat}^*)\\
      0\rar&\mathbf{R}^1\lim\Gamma(X,\sO_X^*)\rar&\HH^1(X,\mathbf{R}\lim\sO_X^*)\rar&\lim \HH^1(X,\sO_X^*)\rar&\mathbf{R}^2\lim\Gamma(X,\sO_X^*).
    \end{tikzcd}
  \end{equation}
  $\bfR^1\lim\sO_X^*$ is the sheafification of the presheaf $U\mapsto\bfR^1\lim\sO_X(U)^*$.  But $X$ has a basis of affinoid perfectoids $U$, on which $\sO_X(U)$ is a perfectoid ring, so $\bfR^1\lim\sO_X(U)^*=0$ (by Proposition \ref{limInts} or \ref{R1LimVanishes} in the integral or Tate cases respectively).  Since the presheaf is 0 on a basis, its sheafification is 0 so that $\Gamma(X,\bfR^1\lim\sO_X^*)=0$ and the top row exhibits an isomorphism:
  \[\HH^1(X,\sO_{X^\tilt}^*)\cong\HH^1(X,\mathbf{R}\lim\sO_X^*).\]
  On the other hand, we know derived limits of abelian groups have cohomology concentrated in degrees 0 and 1 (see for example \cite[Tag 07KW]{stacks-project}), so that the bottom right term of Diagram (\ref{lowTerms}) is 0.  Therefore Sequence (\ref{mainSequence}) emerges as the bottom row.  If we assume further that $\Gamma(X,\sO_X)$ is a perfectoid ring, then the left term of Sequence (\ref{mainSequence}) is 0 (again by Proposition \ref{limInts} or \ref{R1LimVanishes}), completing the proof.
\end{Proof}
\section{Consequences and Examples}\label{consequences}
\subsection{The Untilting Homomorphism}
Let us set some notation.  For an for an abelian group $A$ we denote by $T_p A = \lim A[p^n]$ the $p$-adic Tate module of $A$, and by $A_{\pDiv}$ the maximal $p$-divisible subgroup of $A$, or equivalently, the subgroup of elements with infinite systems of $p^n$'th roots.  In particular, one always obtains an exact sequence $0\to T_pA\to\lim A\to A_{\pDiv}\to 0$.

Now let $X$ be a perfectoid space.  One can compose $\theta:\Pic X^\flat\to\lim\Pic X$ together with the projection $\pi_0:\lim\Pic X\to\Pic X$ onto the first coordinate to obtain the \textit{untilting homomorphism} $\theta_0:\Pic X^\flat\to\Pic X$.  If $\Gamma(X,\sO_X)$ is a perfectoid ring, Theorem \ref{Main} and the previous paragraph tell us that $\theta_0$ fits into following exact sequence:
\[0\longto T_p\Pic X\longto\Pic X^\flat\longtoo{\theta_0}(\Pic X)_{\pDiv}\longto 0.\]
In particular, we have the following corollary to Theorem \ref{Main}.
\begin{Corollary}\label{MainCorollary}
  Let $X$ be a perfectoid space and suppose $\Gamma(X,\sO_X)$ is perfectoid.  Consider the untilting homomorphism $\theta_0:\Pic X^\flat\to\Pic X$.
  \begin{enumerate}[label ={(\arabic*)}]
    \item{$\theta_0$ is injective if and only if $T_p\Pic X=0$.}
    \item{$\theta_0$ is surjective if and only if $\Pic X$ is $p$-divisible.}
    \item{$\theta_0$ is an isomorphism if and only if $\Pic X$ is uniquely $p$-divisible.}
  \end{enumerate}
\end{Corollary}
As a special case of (1), we see that if $\Pic X$ is $p$-torsion free, or more generally if $\Pic X$ has bounded $p^\infty$-torsion, then $\theta_0$ is injective.

As a consequence of the \'etale comparison theorem in prismatic cohomology, Bhatt and Scholze \cite[Corollary 9.5]{bs19} show that a large class of perfectoid affinoid adic spaces have uniquely $p$-divisible Picard groups.  In particular, they show that if $R$ is integral perfectoid, then $R$ and $R[1/p]$ have uniquely $p$-divisible Picard groups.  This gives us the following immediate corollary.
\begin{Corollary}
  Let $R$ be an integral perfectoid ring.  Then
  \begin{eqnarray*}\Pic R^\flat\cong\Pic R&\text{ and }&\Pic (R[1/p]^\flat)\cong\Pic (R[1/p]).\end{eqnarray*}
  In particular, $\theta_0$ is an isomorphism for any perfectoid algebra over a field.
\end{Corollary}
One observes that the same is true in a number of global examples.
\begin{Example}\label{Projectivoid}
  In \cite[Theorem 3.4]{DH21} we computed the Picard group of \textit{projectivoid space}, showing that $\Pic\bP^{n,\perf}\cong\bZ[1/p]$, which is uniquely $p$-divisible.  Therefore the untilting map for projectivoid space is an isomorphism.  Of course, since the value of $\Pic \bP^{n,\perf}$ is independent of perfectoid base field (and therefore tilting) this is unsurprising.
\end{Example}
\begin{Example}\label{Toric}
  Generalizing the previous example, if $X^\perf$ is the perfectoid cover of a smooth proper toric variety $X$ (as in \cite[Section 8]{sh12}), then $\Pic(X^\perf)\cong\Pic(X)[1/p]$, (by \cite[Theorem 4.1]{dhrw20}).  In particular, $p$ acts invertibly on $\Pic(X^\perf)$ and so $\theta_0$ is an isomorphism.
\end{Example}
\begin{Example}\label{AbelianGoodReduction}
  Let $A$ be an abelian variety (over a perfectoid field), and $[p]:A\to A$ the multiplication by $p$ map.  Then passing to the inverse limit there is a perfectoid space $\widetilde{A}\sim\lim A$ (by \cite[Lemme A.16]{ps16} in the case where $A$ has good reduction over an algebraically closed field, and by \cite[Theorem 4.6]{bghsswy18} in general).  Work of Heuer \cite[Theorem 1.10]{heu21b} shows that if $A$ has good reduction, then $p$ acts invertibly in the Picard group of $\widetilde{A}$, so that by Corollary \ref{MainCorollary} we may conclude that $\Pic\widetilde{A}^\flat\cong\Pic\widetilde{A}$.
\end{Example}
We also include an example where $\theta_0$ is neither injective nor surjective, coming from a construction of \cite{bghsswy18} and \cite[Section 6.3]{heu21b} of a perfectoid space whose Picard group is not $p$-divisble nor $p$-torsion free.
\begin{Example}\label{TateCurve}
  Let $K$ be a perfectoid field of characteristic 0, and consider the Tate uniformization of an elliptic curve $E = \bG_m/q^\bZ$, considered as an adic space over $K$.  Suppose $q$ has a coherent system of $p$-power roots$-$for example, if $q$ is in the image of $\sharp:(K^\flat)^*\to K^*$.  Fix such a system so that the symbol $q^{1/p^n}$ makes sense for all $n>0$.  The $p$-power map induces a sequence of isogenies:
  \[\cdots\to\bG_m/q^{1/p^2}\to\bG_m/q^{1/p}\to\bG_m/q,\]
  and by \cite{bghsswy18}, passing to the `tilde' inverse limit produces a perfectoid space:
  \[E_\infty\sim\ilim_n\bG_m/q^{1/p^n}.\]
  By \cite[Proposition 6.1]{heu21b} the Picard group of $E_\infty$ fits into the following exact sequence.
  \[0\longto K^*/q^{\bZ[1/p]}\longto\Pic E_\infty\longto\bZ[1/p]\longto0.\]
  Consider multiplication by $p$ on this sequence:
  \[
  \begin{tikzcd}
    0\rar&K^*/q^{\bZ[1/p]}\rar\ar[d,"\Phi"]&\Pic E_\infty\rar\ar[d,"\Psi"]&\bZ[1/p]\ar[d]\rar&0\\
    0\rar&K^*/q^{\bZ[1/p]}\rar&\Pic E_\infty\rar&\bZ[1/p]\rar&0
  \end{tikzcd}
  \]
  Corollary \ref{MainCorollary} says that the injectivity and surjectivity of the untilting map $\theta_0:\Pic E_\infty^\flat\to\Pic E_\infty$ is controlled by the kernel and cokernel of $\Psi$, and since $\bZ[1/p]$ is uniquely $p$-divisible the snake lemma tells us that these are in turn isomorphic to the kernel and cokernel of $\Phi$.  But these are easy to compute.  Consider the square whose rows are the $p$-power map
  \[
  \begin{tikzcd}
    K^*\rar\dar&K^*\dar\\
    K^*/q^{\bZ[1/p]}\ar[r,"\Phi"]&K^*/q^{\bZ[1/p]}
  \end{tikzcd}
  \]
  The kernel of the composition consists of functions $f\in K^*$ such that $f^p=q^d$ for some $d\in\bZ[1/p]$.  Thus $f/q^{d/p}\in\mu_p$ is a $p$'th root of unity.  In particular, $f$ is given by an element of $q^{\bZ[1/p]}$ and one of $\mu_p$, and furthermore $q^{\bZ[1/p]}\cap\mu_p=1$ so that the kernel is isomorphic to $\mu_p\times q^{\bZ[1/p]}$.  We get the $\ker\Phi$ modding out by the second factor, so that $\ker\Phi\cong\mu_p$.  Arguing for successively higher powers of $p$ we see that the kernel of the projection $\lim\Pic E_\infty\to\Pic E_\infty$ is isomorphic to $\lim_n\mu_{p^n} = \bZ_p(1)$, so that by Theorem \ref{Main}, we have a left exact sequence:
  \begin{equation}\label{tateExact}
    0\longto\bZ_p(1)\longto\Pic E_\infty^\flat\longtoo{\theta_0}\Pic E_\infty,
  \end{equation}
  whose right exactness is controlled by $\coker\Phi$.  The cokernel of the $p$-power map from $K^*$ to itself is $K^*/(K^*)^p$, and we chose $q$ so that $q^{\bZ[1/p]}\subseteq(K^*)^p$ so that:
  \[\coker\Phi = \left(K^*/q^{\bZ[1/p]}\right)/\left((K^*)^p/q^{\bZ[1/p]}\right)\cong K^*/(K^*)^p.\]
  Therefore Sequence \ref{tateExact} is right exact if and only if every element of $K^*$ is a $p$'th power.  In particular, we see that it is right exact if $K$ is algebraically closed, but is not in general$-$for example if $K = \bQ_p\left(p^{1/p^\infty}\right)\hat{ }$.

  In fact, we can give a concrete description of the kernel and cokernel of $\theta_0$ for elements that arise from divisors on the Tate curve.  Observe that the choices of $p$-power roots of $q\in K^*$ determines a unique element $q^\flat\in(K^\flat)^*$.  Then we can compute the tilt of $E_\infty$ as
  \[(E_\infty)^\flat\sim\ilim_n (\bG_m)^\flat/(q^\flat)^{1/p^n}.\]
  Identifying $\colim\Pic E$ as a subset of $\Pic E_\infty$ (and similarly for $\Pic E_\infty^\flat)$, we can identify points of $K^*/q^{\bZ[1/p]}$ (respectively $(K^\flat)^*/(q^\flat)^{\bZ[1/p]}$) with certain degree 0 line bundles on $E_\infty$ (respectively $(E_\infty)^\flat$).  On these line bundles, $\theta_0$ descends from the untilting map $\sharp$.
  \[
  \begin{tikzcd}
    (K^\flat)^*\ar[r,"\sharp"]\dar & K^*\dar\\
    (K^\flat)^*/(q^\flat)^{\bZ[1/p]}\ar[r,"\theta_0"] & K^*/q^{\bZ[1/p]}.
  \end{tikzcd}
  \]
  With this presentation, we see that divisors in the kernel are precisely those with coordinates that map to 1 under $\sharp$ (i.e., elements of $\bZ_p(1)$), elements of the cokernel come from points in $K^*$ that don't have infintely many $p$-power roots.
\end{Example}
\subsection{Trivializing Inverse Systems of Line Bundles in the Analytic Topology}
Let $\sL_0,\sL_1,\sL_2,\cdots$ be a system of line bundles on a perfectoid space $X$ with $\sL_{i+1}^{\otimes p}\cong\sL_i$.  There is an obvious way to construct a pro-\'etale cover of $X$ that simultaneously trivializes all of the $\sL_i$, by further refining \'etale covers trivializing each $\sL_i$ individually and letting $i$ go to infinity.  Nevertheless, it isn't immediately clear that there should be an \'etale cover that simultaneously trivializes all the $\sL_i$.  A consequence of Theorem \ref{Main} is that there is, and in fact there is even an \textit{analytic} cover that does so.
\begin{Corollary}\label{analyticTrivialization}
  Let $X$ be a perfectoid space and $\sL_0,\sL_1,\cdots,$ a system of line bundles on $X$ with $\sL_{i+1}^{\otimes p}\cong\sL_i$.  Then there is an analytic cover $U\to X$ that simultaneously trivializes all the $\sL_i$.
\end{Corollary}
\begin{Proof}
  Although the global sections of $\sO_X$ aren't a priori perfectoid, there is an analytic open cover $V\to X$ where $\Gamma(V,\sO_V)$ is.  Therefore by Theorem \ref{Main}, the inverse system $(\sL_{i}|_{V})\in\lim\Pic V$ is the untilt of a unique $\sL\in\Pic V^\flat$.  There is an analytic cover $U^\flat\to V^\flat$ trivializing $\sL$, which is the tilt of an analytic cover $U\to V$.  By functoriality the inverse system $(\sL_{i}|_{U})\in\lim\Pic U$ is the untilt of $\sL_{U^\flat}$, which is trivial.  Again by Theorem \ref{Main} each $\sL_{i}|_U$ is trivial, completing the proof.
\end{Proof}
\begin{Remark}
  In the case that $X$ is integral or over a field, Corollary \ref{analyticTrivialization} is an easy consequence of \cite[Corollary 9.5]{bs19}.  Indeed, in this case if $U\to X$ is an affinoid cover trivializing $\sL_0$, then their result implies that $\Pic(U)$ is uniquely $p$-divisible.  Therefore, since $\sL_n|_U$ has trivial $p^n$-th power, it must be trivial as well.  Corollary \ref{analyticTrivialization} is slightly more general, allowing for Tate perfectoid spaces that are not over fields (cf. \cite[Example 6.1.5.4]{SW18}).
\end{Remark}
\subsection{$\ell$-adic Cohomological N\'eron-Severi Groups}
We'd like to study how the untilting maps $\theta$ and $\theta_0$ act on N\'eron-Severi groups of perfectoid spaces.  Since we don't have representability of the Picard functor, we need a notion of N\'eron-Severi groups to make sense of this, so we take a cohomological approach.  We begin by defining the $\ell$-adic cycle class map.

Let $X$ be an adic space, and $\ell$ a prime invertible in $\sO_{X}$.  From the long exact sequence on \'etale cohomology associated to the Kummer sequence $1\to\mu_{\ell^n}\to\bG_m\to\bG_m\to1$ one obtains an injection:
\begin{equation}\label{injectionPicModEll}
  \frac{\Pic X}{\ell^n\Pic X}\into\HHe^2(X,\mu_{\ell^n}).
\end{equation}
Passing to the inverse limit among all $n$, one obtains a map $\Pic X\otimes\bZ_\ell(1)\into\HHe^2(X,\bZ_\ell(1))$ and precomposing with the natural map from $\Pic X$, one obtains the \textit{$\ell$-adic cycle class map}:
\[c_{\ell}:\Pic X\to\HHe^2(X,\bZ_\ell(1)).\]
\begin{Definition}
  Let $X$ be an adic space and $\ell$ a prime invertible in $\sO_{X}$.  Then the $\ell$-adic cohomological N\'eron-Severi group is defined to be the image of the $\ell$-adic cycle class map. \[\NS_\ell(X):=\im(c_\ell).\]
\end{Definition}
In certain nice situations, including the case of an abelian variety over an algebraically closed field, this definition agrees with the usual notion of N\'eron-Severi groups.
\begin{Lemma}\label{ladicNSisOK}
  Let $A$ be a proper nonsingular variety over an algebraically closed field of characteristic not equal to $\ell$, and suppose further that $\NS(A)$ is torsion-free.  Then $\NS_\ell(A)\cong\NS(A)$.
\end{Lemma}
\begin{Proof}
  This is well known but we include the proof for completeness.  We first observe that because $\Pic^\circ(A)$ is divisible, we have $\Pic^\circ A\subseteq\ell^n\Pic A$ so that:
  \[\frac{\Pic A}{\ell^n\Pic A}\cong\frac{\Pic A/\Pic^\circ A}{\ell^n\Pic A/\Pic^\circ A}\cong\frac{\NS(A)}{\ell^n\NS(A)}.\]
  In particular, the source of the injection from Equation (\ref{injectionPicModEll}) can be identified with $\NS(A)/\ell^n\NS(A)$ and passing the the inverse limit along $n$ one obtains the following composition whose image agrees with the image of the $\ell$-adic cycle class map $c_\ell$:
  \[\NS(A)\to\NS(A)\otimes\bZ_\ell(1)\into\HHe^2(X,\bZ_\ell(1)).\]
  We finish by observing that the first map is injective because $\NS(A)$ is torsion-free.
\end{Proof}
The $\ell$-adic cohomological N\'eron-Severi group plays well with the untilting homomorphism that is the subject of this paper.
\begin{Proposition}\label{NSThetaInjects}
  Let $X$ be a perfectoid space and $\ell$ prime to the residue characteristic of $X$.  Then the untilting homomorphism $\theta_0$ passes to an injection:
  \[\NS_\ell(X^\flat)\into\NS_\ell(X).\]
\end{Proposition}
\begin{Proof}
  Identify the \'etale sites of $X$ and $X^\flat$.  Since the projection map $\sharp:\bG_m^\flat\cong\lim\bG_m\to\bG_m$ is multiplicative, for all $n$ it induces a map of Kummer sequences
  \[
  \begin{tikzcd}
    1\rar&\mu_{\ell^n}^\flat\rar\ar[d] & \bG_m^\flat\rar\dar & \bG_m^\flat\rar\dar & 1\\
    1\rar & \mu_{\ell^n}\rar & \bG_m\rar & \bG_m\rar & 1.
  \end{tikzcd}
  \]
  Since the $p$-power map on $\mu_{\ell^n}$ is an isomorphism, the projection $\mu_{\ell^n}^\flat\to\mu_{\ell^n}$ is too.  Taking the long exact sequences on \'etale cohomology fits $\theta_0$ into the following diagram.
  \[
  \begin{tikzcd}
    \Pic X^\flat\ar[d,"\theta_0"]\ar[r,"\ell^n"] & \Pic X^\flat\ar[d,"\theta_0"]\ar[r] & \HHe^2(X^\flat,\mu_{\ell^n})\ar[d,"\wr"]\\
    \Pic X\ar[r,"\ell^n"] & \Pic X\rar & \HHe^2(X,\mu_{\ell^n})
  \end{tikzcd}
  \]
  Passing to the limit among all $n$ we obtain the following commutative square.
  \[
  \begin{tikzcd}
    \Pic X^\flat\ar[d,"\theta_0"]\ar[r,"c_{\ell}^\flat"] & \HHe^2(X^\flat,\bZ_\ell(1))\ar[d,"\wr"]\\
    \Pic X\ar[r,"c_\ell"] & \HHe^2(X,\bZ_\ell(1)).
  \end{tikzcd}
  \]
  In particular, the vertical map on the right restricts to an injection between the images of the cycle class maps, completing the proof.
\end{Proof}
Proposition \ref{NSThetaInjects} has the following immediate corollary, arguing analogously to Corollary \ref{MainCorollary}.
\begin{Corollary}
  Let $X$ be a perfectoid space and $\ell$ prime to the residue characteristic of $X$.  If $\Gamma(X,\sO_X)$ is a perfectoid ring, then $T_p\NS_\ell(X)=0$.
\end{Corollary}
\subsection{Appearances of New Line Bundles at Infinite Level}
Both \cite[Theorem 3.4]{DH21} and \cite[Theorem 4.1]{dhrw20} study examples of varieties $X$, together with a `Frobenius like' map $\Phi:X\to X$ so that there is a \textit{mixed characteristic perfection}, a perfectoid space $\widetilde{X}\sim\lim X$ (as the notation suggests, this can be thought of in certain cases as a sort of universal cover).  The content of the theorems in each case is that this limit commutes with taking Picard groups, i.e., that the natural map $\colim\Pic X\to\Pic \widetilde{X}$ is an isomorphism.  Work of Heuer \cite[Section 6.2]{heu21b} shows that this is not true in general, giving an example where this fails in characteristic $p$.  We conclude by giving an example of this failure in characteristic 0.  In both cases, the counterexample consists of an abelian variety $A$ over a perfectoid field and its perfectoid cover $\widetilde{A}\sim\lim_{[p]}A\to A$ (as in Example \ref{AbelianGoodReduction}) such that the induced map $\colim\Pic A\to\Pic \widetilde{A}$ is not surjective.

The idea is to start with an abelian variety whose N\'eron Severi rank jumps modulo $p$, and then use Proposition \ref{NSThetaInjects} to inject the (now larger) N\'eron-Severi group of $\widetilde{A}^\tilt$ into the N\'eron-Severi group of $\widetilde{A}$, thus exhibiting line bundles that cannot come from $A$.  For the perfectoid covers we will need to use the $\ell$-adic cohomological N\'eron-Severi groups introduced above, observing by Lemma \ref{ladicNSisOK} that for abelian varieties over algebraically closed fields the two notions are interchangable.  We first confirm that the N\'eron-Severi rank does not decrease when passing to the perfectoid cover.
\begin{Lemma}\label{phibarinjects}
  Suppose $A$ is an abelian variety over a perfectoid field, and let $\ell$ be prime to the residue characteristic.  Let $\widetilde{A}\to A$ be its perfectoid cover.  The induced map $\colim\NS(A)\to\NS_\ell(\widetilde{A})$ is injective.
\end{Lemma}
\begin{Proof}
  The Kummer sequence $0\longto\mu_{\ell^n}\longto\bG_m\longto\bG_m\longto0$ induces the following diagram whose rows are exact:
  \begin{equation*}\label{diagram}
    \begin{tikzcd}
      \Pic A\dar\ar[r,"\ell^n"] & \Pic A\dar\rar&\HHe^2(A,\mu_{\ell^n})\dar\\
      \colim\Pic A\ar[d]\rar&\colim\Pic A\ar[d]\rar & \colim\HHe^2(A,\mu_{\ell^n})\ar[d,"\eta_n"]\\
      \Pic \widetilde{A}\rar&\Pic \widetilde{A}\rar&\HHe^2(\widetilde{A},\mu_{\ell^n}).
    \end{tikzcd}
  \end{equation*}
  Suppose $(\sL_i)\in\colim\Pic A$ maps to 0 in $\NS_\ell(\widetilde{A})$.  This means that for $n$ large enough, it maps to 0 in $\HHe^2(\widetilde{A},\mu_{\ell^n})$.  By \cite[Corollary 7.8]{sh12}, $\eta_n$ is an isomorphism, so that $(\sL_i)$ maps to 0 in $\colim\HHe^2(A,\mu_{\ell^n})$.  By exactness, $\sL_i$ is an $\ell^n$-th power for $i\gg0$, so that the class of $\sL_i$ in $\NS(A)\cong\bZ^{\rho(A)}$ is a multiple of $\ell^n$ for all $n\gg0$, so that it must be 0.  Therefore $(\sL_i)\in\colim\Pic^\circ A$ and we win.
\end{Proof}
We now have a model for our counterexample.
\begin{Proposition}\label{mainCounter}
  Let $K$ be an algebraically closed perfectoid field with tilt $K^\flat$.  Suppose that $A$ is an abelian variety over $K$, and $B$ an abelian variety over $K^\flat$, whose perfectoid covers satisfy $\widetilde{A}^\flat = \widetilde{B}$.  If $\rho(B)>\rho(A)$, then the map $\colim\Pic A\to\Pic \widetilde{A}$ is not surjective.
\end{Proposition}
\begin{Proof}
  We remind the reader that for any abelian variety $X$ over an algebraically closed field, if $[n]:X\to X$ is multiplication by $n$, then the pullback map $[n]^*:\NS X\to\NS X$ is multiplication by $n^2$ (see for example \cite[2.8(iv)]{mum74}).  Therefore, as $\NS(X)\cong\bZ^{\rho(X)}$, we know $\colim_{[n]^*}\NS(X)\cong\bZ[1/n]^{\rho(X)}$.

  Therefore, the assumptions of the proposition, together with Lemma \ref{phibarinjects} and Proposition \ref{NSThetaInjects} give us the following chain of inequalities:
  \[\rk_{\bZ[1/p]}\colim\NS(A)<\rk_{\bZ[1/p]}\colim\NS(B)\le\rk_{\bZ[1/p]}\NS_\ell(\widetilde{B})\le\rk_{\bZ[1/p]}\NS_\ell(\widetilde{A}).\]
  In particular, $\colim\NS(A)\to\NS_\ell(\widetilde{A})$ cannot surject, and therefore neither can the map in question.
\end{Proof}
We conclude by giving examples of abelian varieties that satisfy the assumptions of Proposition \ref{mainCounter}.  Let $\cA$ be the integral model over $\bZ_p$ of an abelian variety with good reduction over $\bQ_p$.  Let $\bC_p$ be the completion of the algebraic closure of $\bQ_p$, and let $\bC_p^\flat$ be its tilt.  Then one can consider the base change of $\cA$ to $\bC_p$ and $\bC_p^\flat$, which we denote by $A_{\bC_p}$ and $A_{\bC_p^\flat}$ respectively.
\begin{Lemma}\label{ZpTilts}
  In the setup of the previous paragraph, we have
  \[\left(\widetilde{A_{\bC_p}}\right)^\flat\cong \widetilde{A_{\bC^\flat_p}}.\]
\end{Lemma}
\begin{Proof}
  Let $\varpi^\flat$ be a pseudouniformizer for $\bC_p^\flat$ and $\varpi = (\varpi^\flat)^\sharp$.  By the tilting equivalence \cite{sh12}, both $\widetilde{A_{\bC_p}}$ and $\widetilde{A_{\bC_p^\flat}}$ are determined (up to almost isomorphism) by their models over $\cO_{\bC_p^\flat}/\varpi^\flat = \cO_{\bC_p}/\varpi = :R$, which extends $\bF_p$.  In each case, we observe that this model must be the scheme $\ilim_{[p]}\left(A_{\bF_p}\times_{\spec\bF_p}\spec R\right)$.
\end{Proof}
To construct abelian varieties satisfying the assumptions of Proposition \ref{mainCounter}, we may therefore start with abelian varieties with good reduction over $\bQ_p$.  Abelian varieties over $\bQ_p$ whose N\'eron-Severi ranks increase upon reduction modulo $p$ are abundant.  Take, for example, $A = E\times E$ where $E$ is a non-CM elliptic curve over $\bQ_p$ with supersingular reduction.  In this case we have $\rho(A_{\bC_p^\flat}) = 6>3 = \rho(A_{\bC_p})$, giving the desired example.
\bibliography{bib}{}
\bibliographystyle{alpha}
\end{document}